# Cycles and Patterns in the Sieve of Eratosthenes

George F Grob and Matthias Schmitt

Revised October 25, 2019

We describe recurring patterns of numbers that survive each wave of the Sieve of Eratosthenes, including symmetries, uniform subdivisions, and quantifiable, predictive cycles that characterize their distribution across the number line. We generalize these results to numbers that are relatively prime to arbitrary sets of prime numbers and derive additional insights about the distribution of integers counted by Euler's $\phi$-function.

## 1. Introduction

In the third century B.C., Eratosthenes of Cyrene (276– 194 BC) invented his famous sieve to simplify the identification of prime numbers. It has remained the backbone of prime number counting and analysis ever since.

In this article, we take a different perspective on the sieve. Rather than using it as a source of prime numbers and counting them, we view its iterative waves as intriguing worlds unto themselves and describe the distribution and relationships of the numbers that remain within them. We call the survivors of its $n^{th}$ iteration of the sieve *n-primes*.

Much is already known about n-primes. One of us wrote about them some fifty years ago[i], and the other did so more recently, calling them "relative primes."[ii] Many others, too numerous to list here, have also studied them, calling them "prime modular numbers," "wheel factorization spokes[iii]" or other terms. By these and other names, mathematicians of the past and still today use n-primes as starting points for sophisticated theorems about the distribution of primes, including twin primes.

For readers new to these concepts or needing a quick refresher, we start with definitions of and well known fundamental relationships among n-primes, identifying and describing overarching structures that occur within each succeeding wave of Eratosthenes's sieve. This is the subject of theorems 1 – 3.

Theorem 4, about twin n-primes, has also been independently discovered or more fully developed by others. It is important to note that the twin prime conjecture (regarding the likely infinite number of pairs of primes separated by 2, is *not* a focal point of this article. Rather, theorem 4 is included, as it was in the original article, as a reminder of a notable feature of the distribution of n-primes. That said, we provide in corollary 2 a generalization about twin n-primes, broadening the analysis to include prime pairs separated by integers greater than 2.

Theorem 5 is about transitions from one wave of Eratosthenes's sieve to the next.

In theorems 6 – 9 we provide new theorems that identify and hopefully provide insights about subordinate cycles within the waves of Eratosthenes's sieve. From all these theorems we generalize our results to arbitrary sets of primes (including infinite sets), relatively prime factors of composite numbers, and patterns of integers counted by Euler's $\phi$-function.

To be clear, we distinguish between the well-known "waves" of Eratosthenes's sieve itself (i.e., the survivors of the first, second, third, and so on, through the $n^{th}$ iteration of the sieve) and subordinate cycles of integers which appear within each wave. The waves are widely known among mathematicians today. *It is the subordinate cycles of n-primes within the overarching waves of Eratosthenes's sieve, beginning with theorem 6, and generalizations of well-known theorems about the sieve to arbitrary sets of relatively prime integers and Euler's $\phi$-function* that are the main focus of this paper.

**2. Background: N-Primes and Their Waves**

Definition (n-primes): Label the n-th prime number as $p_n$. Thus, $p_1 = 2$, $p_2 = 3$, $p_3 = 5$, $p_4 = 7$, $p_5 = 11$, etc. If an integer $a$ is not divisible by any of the first n primes, let us call $a$ an n-prime. In essence, then, an n-prime is a survivor of the first n iterations of Eratosthenes's sieve.

Define $\Pi p_i = p_1 \cdot p_2 \cdot p_3 \ldots p_n$, the product of the first n primes. If the range of i is other than 1 ... n, or if needed or useful for clarity or emphasis, we will specify it appropriately.

Let $Z_0$ be the set of non-negative natural numbers 0, 1, 2 . . . etc.

Let $P_n = \{x \text{ in } Z \text{ such that } 0 \leq x < \Pi p_i, 1 \leq i \leq n\}$ be the set of integers smaller than $\Pi p_i$.

Let $\phi(x)$ = Euler's $\phi$ function

**Theorem 1.** An integer $a$ is an n-prime if and only if $K \cdot \Pi p_i + a$ is an n-prime $_{1 \leq i \leq n}$

for all $K \in Z$.

Proof. Suppose $a$ is an n-prime but that $K \cdot \Pi p_i + a$ is not an n-prime for some $K \in Z$. Then $p_o | (K \cdot \Pi p_i + a)$ for some $p_1 \leq p_o \leq p_n$. But $p_o | K \cdot \Pi p_i$. So, $p_o | (K \cdot \Pi p_i + a - K \cdot \Pi p_i)$. Thus, $p_o | a$, a contradiction to the fact that $a$ is an n-prime. The converse is proved in the same way.

Thus, once we know the distribution of n-primes in $P_n$, then, by Theorem 1, the pattern established there forms a cycle that will be repeated every $\Pi p_i$ consecutive integers.

We now investigate the distribution of n-primes within $P_n$.

**Theorem 2.** An integer $a$ in $P_n$ is an n-prime if and only if $\Pi p_i - a$ is an n-prime. That is, the n-primes are arranged symmetrically within $P_n$.

Proof. The proof is similar to that for Theorem 1. Suppose $a$ is an n-prime but that $\Pi p_i - a$ is not an n-prime for some $a$. Then

$p_i | \Pi p_i - a$, for some $i = 1, 2 \ldots n$.

But $p_i | \Pi p_i$. So, $p_i | \Pi p_i - (\Pi p_i - a)$. Thus, $p_i | a$, a contradiction to the fact that $a$ is an n-prime.

**Theorem 3.** There are $(p_2 - 1)(p_3 - 1) \ldots (p_n - 1)$ n-primes in $P_n$.

Proof. The number of n-primes in $P_n$ is equal to the number of positive integers relatively prime to and smaller than $\Pi p_i$. Using Euler's $\phi$ function, there are

$\phi(\Pi p_i) = (p_1 - 1)(p_2 - 1)(p_3 - 1) \ldots (p_n - 1)$ such integers.

Since $p_1 = 2$, and $2 - 1 = 1$, this is reduced to

$\phi(\Pi p_i) = (p_2 - 1)(p_3 - 1) \ldots (p_n - 1)$.

Again, by theorem 1, each of the infinitely repeating cycles of $\Pi p_i$ integers contain that same number of n-primes, and in the same symmetrical order as the first one.

A key tool for our analysis throughout this paper will be a function

$f_n(x) =$ the number of n-primes $\leq x$.

We will have more to say about this function in subsequent sections. For now, we use it to summarize theorems 1 – 3 as follows.

**Corollary 1.** For any integer $K \geq 0$,

$f_n(K \cdot \Pi p_i \pm x) = f_n(K \cdot \Pi p_i) \pm f_n(x) = K \cdot (p_2 - 1)(p_3 - 1) \ldots (p_n - 1) \pm f_n(x)$

**Examples.** The following are examples of the application of theorems 1 – 3 and Corollary 1, in the case of n = 4, i.e., the 4-primes (non-negative integers that are not divisible by 2, 3, 5, and 7). By theorem 1, their distribution within the integers is repeated every $\Pi p_i = 210$ times.
$\scriptstyle 1 \leq i \leq 4$

As predicted by theorem 2, they are symmetrically distributed within this interval (e.g., 1 and 209; 31 and 179; 103 and 107).

As predicted by theorem 3, there are $(2 - 1)(3 - 1)(5 - 1)(7 - 1) = 1 \cdot 2 \cdot 4 \cdot 6 = 48$ of them.

Table 1 lists the 4-primes within this first cycle of 210 integers.

**Table 1**
**The first cycle of 4-primes**

| 1 | 11 | 13 | 17 | 19 | 23 |
|---|----|----|----|----|----|
| 29 | 31 | 37 | 41 | 43 | 47 |
| 53 | 59 | 61 | 67 | 71 | 73 |
| 79 | 83 | 89 | 97 | 101 | 103 |
| 107 | 107 | 113 | 121 | 127 | 131 |
| 139 | 139 | 143 | 149 | 151 | 157 |
| 163 | 167 | 169 | 173 | 179 | 181 |
| 187 | 191 | 193 | 197 | 199 | 209 |

## 3. Twin N-Primes

Definition (twin n-primes): A famous unanswered question about primes is: How many twin primes are there? By twin primes we mean a pair of primes $a$ and $a + 2$. We can similarly define a pair of twin n-primes to be a pair of n-primes $a$ and $a + 2$. The following theorem addresses the question of how many pairs of twin n-primes are found within each cycle.

**Theorem 4.** There are

$$\prod_{3 \leq i \leq n}(p_i - 2) = (p_3 - 2)(p_4 - 2) \cdots (p_n - 2)$$

integers $x$ in $P_n$ such that $x - 1, x + 1$ is a pair of twin n-primes.

Proof. We base our proof on cross product multiplication of the modular arithmetic sets for the first n primes. Let $A_i = (x$ in Z such that $0 \leq x \leq p_i - 1)$. Define f: mapping Z onto the cross product $A_1 \times A_2 \times \ldots \times A_n$ by for $a$ in Z, $f(a) = (a_1, a_2, \ldots, a_n)$ where for $1 \leq i \leq n$, $a_i$ is the remainder when $a$ is divided by $p_i$. By the Chinese remainder theorem, f is onto. Given the n-tuple $a' = (a_1, a_2 \ldots a_n)$, let $a$ be such that $f(a) = a'$. Then, since $p_i$ divides $K \cdot \Pi p_i$ for all $1 \leq i \leq n$, $f(K \cdot \Pi p_i + a) = a'$ for all K in Z. Choose $K_o$ such that $1 \leq f(K_o \cdot \Pi p_i + a) \leq \Pi p_i$. Then for each $a'$ in $(A_1 \times A_2 \times \ldots \times A_n)$, there is some $1 \leq a \leq \Pi p_i$ such that $f(a) = a'$. Let f* be the restriction of f to $P_n$. Then f* mapping $P_n$ onto $A_1 \times A_2 \times \ldots \times A_n$ is onto; and since card $(A_1 \times A_2 \times \ldots \times A_n) = \Pi p_i = $ card$(P_n)$, the mapping, f* is a bijection.

It is clear that $x \in Z$ is an n-prime if and only if $f(x)$ has no zero entries. Then for $x$ in $P_n$, the pair $x - 1, x + 1$ is a pair of twin n-primes if and only if $f^*(x) = f(x)$ never has 1 or $p_i - 1$ as its i-th entry, where $1 \leq i \leq n$. By counting all n-tuples without such entries, the theorem is proved.

**Corollary 2.** More generally, for any integer $a \neq 0$ (mod $p_i$, $1 \leq i \leq n$) there are
$$\prod_{3 \leq i \leq n}(p_i - 2) = (p_3 - 2)(p_4 - 2) \cdots (p_n - 2)$$

integers $x$ in $P_n$ such that the pair $x - a, x + a$ (mod $P_n$), are both n-primes. This is a simple generalization of theorem 4. Interestingly, it can be generalized even further, outside the boundaries of modular arithmetic, to any integer, a, no matter how large, if we expand our consideration of the domain of x to any x > a; or if we extend the domain of our modular arithmetic across the entire number line, including negative integers. And again, any of these patterns will repeat every $K \cdot \Pi(p_i)$ intervals for every integer $K \geq 1$.

Corollary 2 may be further generalized, holding true for any two integers a and b that are not congruent mod pi, for all i, $1 \leq i \leq n$, That is, for any two such integers, a and b, (even if one of them is congruent to 0 (mod pi, for any i, $1 \leq i \leq n$) there are
$$\prod_{3 \leq i \leq n}(p_i - 2) = (p_3 - 2)(p_4 - 2) \cdots (p_n - 2)$$

integers $x$ in $P_n$ such that the pair $x - a, x + b$ (mod $P_n$), are both n-primes.

Should a be congruent to b for any $p_i$, then the formula holds true by substituting $(p_i - 1)$ for $(p_i - 2)$ for any such $p_i$.

In the next section we provide examples of the distribution of n-primes, twin n-primes, and their repetitive cycles.

See examples of theorem 4, twin 4-primes, in table 1, (e.g., 11 and 13, 17 and 19, 71, and 73). Similarly, see examples of corollary 2, with a = 2, (e.g., 19 and 23, 67 and 71, and 139 and 143); and with a = 3 (e.g., 31 and 37, 47 and 53, and 61 and 67).

**4. Transitions from One Wave to the Next**

We may think of Eratosthenes's sieve as a kind of killing machine, with the first survivor of the last "kill" as the generator of the next round, wiping out an infinite number of survivors of the previous round, but still leaving infinitely more survivors, albeit somewhat "thinned out" from the previous killing.

As noted in Theorems 1 – 3, an overarching feature of such surviving integers is a pattern that repeats itself every $\Pi p_i = p_1 \cdot p_2 \cdot p_3 \cdots p_n$ integers. So, to understand the effect of each iteration of Eratosthenes's sieve, it is useful to describe the pattern of the n-primes within the first $\Pi p_i$ numbers, as well as how those patterns change with each iteration of the sieve. Of course, $\Pi p_i$ grows factorially, actually much faster, with n. Thus, the patterns of growth and distribution are difficult to visualize.

To overcome this problem, we use a combination of rigorous algebraic proofs that apply to all "n" while simultaneously illustrating the results through theorems and graphical representations of what happens between first 3 n-primes ($p_1 = 2$, $p_2 = 3$, $p_3 = 5$) and the 4$^{th}$ one, $p_4 = 7$. To begin, consider the following graphical presentation of the 3-primes, repeated 7 times.

**Figure 1**
**The First 7 Cycles of 3-Pimes**
**(Integers not divisible by 2, 3, or 5)**

| 30 | 60 | 90 | 120 | 150 | 180 | 210 |
|---|---|---|---|---|---|---|
| 1  | 31 | 61 | 91  | 121 | 151 | 181 |
| 7  | 37 | 67 | 97  | 127 | 157 | 187 |
| 11 | 41 | 71 | 101 | 131 | 161 | 191 |
| 13 | 43 | 73 | 103 | 133 | 163 | 193 |
| 17 | 47 | 77 | 107 | 137 | 167 | 197 |
| 19 | 49 | 79 | 109 | 139 | 169 | 199 |
| 23 | 53 | 83 | 113 | 143 | 173 | 203 |
| 29 | 59 | 89 | 119 | 149 | 179 | 209 |

The columns represent increments of 30 (= $P_3$) with integers relatively prime to the first 3 primes, 2, 3, and 5. The columns are repeated 7 ($p_4$) times, covering the domain, 210 (= $P_4$), of the first 4 primes, but before 7 wipes out its multiples.

There are 8 3-primes within each interval of 30, i.e., $P_3$. In each subsequent interval of 30, represented by the remaining six columns, the 3-primes are distributed exactly the same, but with $\Pi p_i$, (i = 1⋯3) = 30, added in each bar.

The first thing that bears repeating here is that within $P_n$, the set of the first $\Pi p_i$ numbers, the n-primes are located symmetrically. This is illustrated in each column, where the first four 3-primes are mirror images of the last four, and the last four are mirror images of the first four in the next column, and so on. This graphically illustrates theorems 1 - 4 for the case of n = 3.

> Theorems 1 & 2: An integer *a* is an n-prime if and only if $K \cdot \Pi p_i + a$, and $K \cdot \Pi p_i - a$ $(1 \leq i \leq n)$, are n-primes for all $K \in Z$.
> Theorem 3: There are $(p_2 - 1)(p_3 - 1) \ldots (p_n - 1)$ n-primes in $P_n$, i.e., in this case, there are $(2 - 1)(3 - 1)(5 - 1) = 1 \cdot 2 \cdot 4 = 8$ 3-primes < 30.
> Theorem 4: There are $(p_3 - 2) = (5 - 2) = 3$ integers *a* in $P_3$ such that $x - 1, x + 1$ is a pair of twin 3-primes, namely 12, 18, and 30.

We now move to the more complex question of what happens to these 3-primes when we eliminate those that are divisible by the next prime, $p_4$ (i.e., 7). The result is a list of 4-primes within its domain of interest, namely the first 210 integers, which is equal to $P_4$, or $2 \cdot 3 \cdot 5 \cdot 7 = 210$. ( As noted earlier, the resulting pattern will be repeated with every multiple of 210.)

**Figure 2**
**3-Pimes deleted by $p_4 = 7$**
**Leaving integers not divisible by 2, 3, 5, or 7 = 4-primes**

| 30 | 60 | 90 | 120 | 150 | 180 | 210 |
|---|---|---|---|---|---|---|
| 1  | 31 | 61 | **X** | 121 | 151 | 181 |
| **X** | 37 | 67 | 97 | 127 | 157 | 187 |
| 11 | 41 | 71 | 101 | 131 | **X** | 191 |
| 13 | 43 | 73 | 103 | **X** | 163 | 193 |
| 17 | 47 | **X** | 107 | 137 | 167 | 197 |
| 19 | **X** | 79 | 109 | 139 | 169 | 199 |
| 23 | 53 | 83 | 113 | 143 | 173 | **X** |
| 29 | 59 | 89 | **X** | 149 | 179 | 209 |

In each column, an X marks the spot where a 3-prime is wiped out by a multiple of the 4$^{th}$ prime, 7. The number of surviving 3-primes in $P_4$ is $(2 - 1)(3 - 1)(5 - 1)(7 - 1) = 1 \cdot 2 \cdot 4 \cdot 6 = 48$, as Theorem 3 predicts.

Similarly, as Theorem 4 predicts, there are $(p_3 - 2)(p_4 - 2) = (5 - 2)(7 - 2) = 3 \cdot 5 = 15$ integers *x* in $P_4$ such that $x - 1, x + 1$ is a pair of twin 4-primes, namely 12, 18, 30, 42, 60, 72, 102, 108, 138, 150, 168, 180, 192, 198, and 210.

Furthermore, just as was the case for $\Pi p_i$, $(1 \leq i \leq 3)$, the distribution of the remaining 4-primes (as well as the deleted 3-primes) are arranged symmetrically with the interval $\Pi p_i = 210$, $(1 \leq i \leq 4)$, as our theory predicts will happen for any n.

It is also worth noting that in the example above that the 4$^{th}$ prime, 7, knocks out exactly one -prime in each of the 8 rows of the diagram. This is a generalizable pattern as proven in Theorem 5 below.

**Theorem 5.** The $(n+1)^{th}$ prime, $p_{n+1}$, divides one and only one n-prime in the series $K \cdot \Pi p_i + a$ for any n-prime *a* and $0 \leq K \leq p_{n+1} - 1$.
$_{1 \leq i \leq n}$

Proof. The proof is by contradiction.

Suppose $p_{n+1}$ divides both $a \cdot \Pi p_i + x$ and $b \cdot \Pi p_i + x$, where a and b are integers between 1 and $p_{n+1}$ and $a > b$. Then $p_{n+1}$ divides $(a \cdot \Pi p_i + x) - (b \cdot \Pi p_i + x) = (a-b) \cdot \Pi p_i$. However, $p_{n+1} > a-b$ and so cannot divide it. Similarly, $p_{n+1} >$ any of the prime factors of $\Pi p_i$, which are all primes $\leq p_n$, and so cannot divide it. Thus $p_{n+1}$ cannot divide more than 1 n-prime of the form $K \cdot \Pi p_i + x$ for any integer.

On the other hand, $p_{n+1}$ must divide at least one of the series $K \cdot \Pi p_i + x$ where $K = 0 \cdots p_{(n+1)}-1$ and x is a given n-prime with $0 < x \leq \Pi p_i$. If not, for one or more K, all the n-primes in those rows would also be (n+1)-primes, and thus for those rows, there would be n+1 of them. Based on our earlier findings, there are $(2-1) \cdot (3-1) \cdots (p_n-1)$ n-primes in each column, and noting that there are n+1 columns there would be more than $(2-1) \cdot (3-1) \cdots (p_{n+1}-1)$
(n+1)-primes altogether, which contradicts our earlier theorem about the number of (n+1) primes $\leq \Pi p_i$, $1 \leq i \leq n+1$.

**5. Cycles Within the Waves**

The previous theorems reveal the overarching patterns of n-primes that remain after each wave of Eratosthenes's sieve—their repetitive, symmetrical patterns of predictable size; the existence and predictable numbers of closely related "twin pairs"; and preservation of many but not all distribution patterns of the previous wave. However, the sizes of each wave of n-primes rise so dramatically with each "n", that distribution patterns cannot be easily distinguished visually. This is illustrated in Table 1.

**Table 1**
**Size of domains of interest for the first 10 n-primes**

|     | Primes | Size of Each Wave $\Pi p_i$ | Number of n-primes in Each Wave $\Pi (p_i - 1)$ |
|-----|--------|-----------------------------|--------------------------------------------------|
| 1.  | 2      | 2                           | 1                                                |
| 2.  | 3      | 6                           | 2                                                |
| 3.  | 5      | 30                          | 8                                                |
| 4.  | 7      | 210                         | 48                                               |
| 5.  | 11     | 2,310                       | 480                                              |
| 6.  | 13     | 30,030                      | 5,760                                            |
| 7.  | 17     | 510,510                     | 92,160                                           |
| 8.  | 19     | 9,699,690                   | 1,658,880                                        |
| 9.  | 23     | 223,092,870                 | 36,495,360                                       |
| 10. | 29     | 6,469,693,230               | 1,021,870,080                                    |

Given the enormity of these and subsequent waves, we looked for patterns within the waves in order to more intuitively and mathematically describe their interiors. We found such patterns within each wave of Eratosthenes's sieve, that, while still ultimately quite large, are easily predictable, calculatable, and conceptually revealing. We call these the "cycles" within the "waves."

A key to our analyses of these cycles is a formula first formulated by Ernst Meissel in 1870, namely

(1) $f_n(x) = f_{n-1}(x) - f_{n-1}(x/p_n)$

where $f_n(x)$ is the number of n-primes $\leq x$. The basis of this formula is straightforward – the $n^{th}$ iteration of Eratosthenes's sieve picks off those multiplicands of $p_n$ that are $\leq x/p_n$ of those integers that are remaining from the previous, $(n-1)^{st}$, iteration of the sieve.

We use Meissel's formula as a basis for many of our findings and proofs, including the following theorems which describe well defined and recurring cycles n-primes within each successive wave of Eratosthenes's sieve.

It is worth noting that while the result of Meissel's formula is always an integer (a number of primes), the domain of its argument is the set of non-negative real numbers, R, i.e.,

$$f_n: R \longrightarrow Z$$

This is obvious in its last term, $f_{n-1}(x/p_n)$, which, even if x were an integer, $x/p_n$ would not be unless $p_n \mid x$. But the formula makes sense and is useful even more generally in its application to real numbers x, especially where x is real number forming the boundary between two domains of n-primes, which is the next topic of our analysis.

Each wave of Eratosthenes's sieve is divided into cycles of equal length, each containing the same numbers of n-primes. More specifically, each $P_n$ may be uniformly divided into cycles length $\prod_{1 \leq i \leq n} p_i /(p_j - 1)$ for each prime $p_j$, $1 \leq j \leq n$, each with $f_{n-1}(\prod_{1 \leq i \leq n-1} p_i)$ n-primes.

We prove this in four stages, through theorems 6 – 9.

**Theorem 6.**  $f_n((\prod_{1 \leq i \leq n} p_i)/(p_n - 1)) = f_{n-1}(\prod_{1 \leq i \leq n-1} p_i)$ for any n.

For example,

$$(2)\ f_4((\prod_{1 \leq i \leq 4} p_i)/(p_4 - 1))$$
$$= f_4(2 \cdot 3 \cdot 5 \cdot 7/(7-1))$$
$$= f_4(210/6)$$
$$= f_4(35)$$
$$= 8$$
$$= f_3(2 \cdot 3 \cdot 5)$$
$$= f_3(30)$$
$$= f_3(\prod_{1 \leq i \leq 3} p_i)$$

Proof.

$$(3)\ f_n((\prod_{1 \leq i \leq n} p_i)/(p_n - 1))$$
$$= f_{n-1}((\prod_{1 \leq i \leq n} p_i)/(p_n - 1)) - f_{n-1}((\prod_{1 \leq i \leq n} p_i)/(p_n - 1)p_n)\ \text{[from (1)]}$$
$$= f_{n-1}((\prod_{1 \leq i \leq n-1} p_i) \cdot p_n /(p_n - 1)) - f_{n-1}((\prod_{1 \leq i \leq n-1} p_i) \cdot p_n /(p_n - 1) \cdot p_n)$$

[switching $(\prod_{1 \leq i \leq n-1} p_i) \cdot p_n$ for $(\prod_{1 \leq i \leq n} p_i)$]

$$= f_{n-1}((\prod_{1 \leq i \leq n-1} p_i) \cdot (1 + 1/(p_n - 1)) - f_{n-1}((\prod_{1 \leq i \leq n-1} p_i)/(p_n - 1))\ \text{[switching 1 +}$$

$1/(p_n - 1)$ for $p_n /(p_n - 1)$ and canceling $p_n$ in both numerator and

denominator of the second term]

$= f_{n-1}(\prod_{1 \leq i \leq n-1} p_i) + f_{n-1}((\prod_{1 \leq i \leq n-1} p_i)/(p_n - 1)) - f_{n-1}((\prod_{1 \leq i \leq n-1} p_i)/(p_n - 1))$

[from corollary (1) with "n-1" replacing "n"]

$= f_{n-1}(\prod_{1 \leq i \leq n-1} p_i)$ [dropping the last two mutually cancelling terms]

We may generalize theorem 6 to integral multiples K of $(P_n/(p_n - 1))$. That is,

**Theorem 7.** For any integer $K \geq 1$, $f_n(K \cdot (\prod_{1 \leq i \leq n} p_i)/(p_n - 1)) = f_{n-1}(K \cdot (\prod_{1 \leq i \leq n-1} p_i))$.

Proof. The proof is achieved by introducing $K \cdot (\prod_{1 \leq i \leq n} p_i)$ and $K \cdot (\prod_{1 \leq i \leq n-1} p_i)$ in proof (3).

For example, let K = 2, introducing it into (2), as follows:

(4)  $f_4(2 \cdot \prod_{1 \leq i \leq 4} p_i) / (p_4 - 1))$

$= f_4(2 \cdot 2 \cdot 3 \cdot 5 \cdot 7 / (7-1))$
$= f_4(2 \cdot 35)$
$= f_4(70)$
$= 16$
$= f_3(2 \cdot \prod_{1 \leq i \leq 3} p_i)$

Continuing in this fashion, we find the following sequence for K = 1⋯ 6.

(5)  $f_4(1 \cdot 210)/6) = f_4(35) = 8$
$f_4(2 \cdot 210)/6) = f_4(70) = 16$
$f_4(3 \cdot 210)/6) = f_4(105) = 24$
$f_4(4 \cdot 210)/6) = f_4(140) = 32$
$f_4(5 \cdot 210)/6) = f_4(175) = 40$
$f_4(6 \cdot 210)/6) = f_4(210) = 48$

Thus, we see that the 4-primes (numbers relatively prime to 2, 3, 5, and 7) are evenly divided within the domain of 4-primes, namely 210, into six intervals of 35, each with 8 4-primes. This is illustrated in figure 3.

**Figure 3**
**Integers not divisible by 2, 3, 5, or 7**
**6 Equal Intervals of 35, each with 8 4-primes**

| 35 | 70 | 105 | 140 | 175 | 210 |
|---|---|---|---|---|---|
| 1  | 37 | 71  | 107 | 143 | 179 |
| 11 | 41 | 73  | 109 | 149 | 181 |
| 13 | 43 | 79  | 113 | 151 | 187 |
| 17 | 47 | 83  | 121 | 157 | 191 |
| 19 | 53 | 89  | 127 | 163 | 193 |
| 23 | 59 | 97  | 131 | 167 | 197 |
| 29 | 61 | 101 | 137 | 169 | 199 |
| 31 | 67 | 103 | 139 | 173 | 209 |

**Additional Uniform Subdivisions of n-primes Using a Generalization of Meissel's Formula**

To recap, we have just shown that we can divide any set of integers

$P_n = \{x \text{ in } Z \text{ such that } 0 \leq x < \prod_{1 \leq i \leq n} p_i\}$ into $(p_n - 1)$ equal segments, each one of which contains the same number of n-primes, namely $f_{n-1}(\prod_{1 \leq i \leq n-1} p_i) = \prod_{1 \leq i \leq n-1}(p_i - 1)$.

This raises the question of whether there are other segments of $P_n$, each of equal length, and each containing the same number of n-primes. The answer is yes, but to understand how this is possible, we need to revisit the sieve of Eratosthenes and the numbers that survive the $n^{th}$ iteration of it, numbers we are calling "n-primes."

It turns out that if we stop with the $n^{th}$ iteration of the sieve, the results do not depend on the order in which the sieve was applied. Eratosthenes started with the first prime, 2. Then he used the smallest survivor of this first round, namely 3, as the agent of killing off the next generation of composite numbers that are multiples of it. Then, he used the smallest survivor of that second round of destruction, namely 5, as the agent of the third generation, and so on.

That said, as far as the final result of the $n^{th}$ killing is concerned, the results would have been the same no matter what order of the first n primes he used to knock them off. So, for example, he could have first wiped out all multiples of 5. Of those that survived that round, he could have eliminated remaining multiples of 2, then of, say, 7, and finally of 3. Any sequence of killing off survivors of the first 4 primes would give the same results—i.e. the survivors would be the same as what we are calling the 4-primes, i.e., integers not divisible by any of the first four primes (2,3,5, and 7).

Based on this, we can generalize Meissel's formula, $f_n(x) = f_{n-1}(x) - f_{n-1}(x/p_n)$, as follows.

**Theorem 8.** $f_{Mn}(x) = f_{Mn-pi}(x) - f_{Mn-pi}(x/p_i)$ for any $p_i \,\varepsilon\, M_n$

Proof. Let $M_n$ be the set of the first n primes, and let $p_i$ be any prime in $M_n$. And so, let $M_n-p_i$ be the set of all the primes in $M_n$ other than $p_i$.

Similarly, define an $M_n$-prime as an integer relatively prime to all the primes in $M_n$. (Thus, an $M_n$-prime is the same thing as an n-prime). Similarly, define an $M_n-p_i$ prime as an integer relatively prime to all the primes in $M_n-p_i$. And thus, define

(6) $f_{Mn}(x)$ = the number integers $\leq x$ that are not divisible by any primes in $M_n$. (Thus, $f_{Mn}(x) = f_n(x)$) and

(7) $f_{Mn-pi}(x)$ = the number integers $\leq x$ that are not divisible by any primes in $M_n-p_i$

As in the original formulation of Meissel's formula, each $i^{th}$ iteration of Eratosthenes's sieve picks off those multiplicands of $p_i$ that are $\leq x/p_i$ of those integers that are remaining from the previous iteration of the sieve. Thus, the proof flows directly from the definitions and the sequence of the application of the sieve.

We return now to our question of whether there are other segments of $P_n$, each segment of equal length, and each of which contains the same number of n-primes. The answer is yes, based on applying the generalization (8) to the case of $M_n = P_n$. Using the same arguments and methods that have brought us to this point describing n-primes and associated formulas for counting them, we may conclude, for example, choosing $p_i = p_3 = 5$,

(8) $f_{Mn} ((\Pi p_i)/(p_3 - 1))$
$\quad\quad\quad {}_{1 \leq i \leq 4}$

$\quad = f_{Mn} (2 \cdot 3 \cdot 5 \cdot 7 / (5-1))$
$\quad = f_{Mn} (210/4)$
$\quad = f_{Mn} (52.5)$ (i.e., the number of integers $\leq 52.5$ relatively prime to 2, 3, 5, and 7; namely, 1, 11, 13, 17, 19, 23, 29, 31, 37, 41, 43, and 47)
$\quad = 12$
$\quad = f_{Mn-p3} ((\Pi p_i)/p_3))$
$\quad\quad\quad\quad {}_{1 \leq i \leq 4}$
$\quad = f_{Mn-p3} (2 \cdot 3 \cdot 7)$
$\quad = f_{Mn-p3} (42)$ (i.e., the number of integers $\leq 42$ relatively prime to 2, 3, and 7; namely, 1, 5, 11, 13, 17, 19, 23, 25, 29, 31, 37, and 41)

Thus, there are 12 4-primes within the first quarter of $P_4$, that are $\leq 52.5$.

As before, we may generalize Theorem 7 to integral multiples K of $((\Pi p_i)/(p_n - 1))$, to

**Theorem 9.** $f_{Mn} (K(\Pi p_i)/(p_j - 1)) = f_{Mn-pj} (K(\Pi p_i))$
$\quad\quad\quad\quad\quad {}_{1 \leq i \leq n} \quad\quad\quad\quad\quad {}_{i \in Mn-pj}$

For example, we find the following sequence for $n = 4$, $p_j = 5$, $K = 1 \cdots 4$.

(9) $f_4 (1 \cdot 210)/4) = f_4 (52.5) = 12$
$\quad f_4 (2 \cdot 210)/4) = f_4 (105) \phantom{.0}= 24$
$\quad f_4 (3 \cdot 210)/4) = f_4 (157.5) = 36$
$\quad f_4 (4 \cdot 210)/4) = f_4 (210) \phantom{.0}= 48$

We illustrate this in Figure 4.

**Figure 4**
**Integers not divisible by 2, 3, 5, or 7**
**4 Equal Intervals of 52.5 (= $2 \cdot 3 \cdot 5 \cdot 7 / (5 - 1)$), each with 12 4-primes**

| 52.5 | 105 | 157.5 | 210 |
|---|---|---|---|
| 1 | 53 | 107 | 163 |
| 11 | 59 | 109 | 167 |
| 13 | 61 | 113 | 169 |
| 17 | 67 | 121 | 173 |
| 19 | 71 | 127 | 179 |
| 23 | 73 | 131 | 181 |
| 29 | 79 | 137 | 187 |
| 31 | 83 | 139 | 191 |
| 37 | 89 | 143 | 193 |
| 41 | 97 | 149 | 197 |
| 43 | 101 | 151 | 199 |
| 47 | 103 | 157 | 209 |

The effect of Theorem 9 is to verify what we foretold in the introduction to Theorems 6 – 9, namely that each wave of the Eratosthenes's sieve is divided into cycles of equal length, each containing the same numbers of n-primes. Specifically, $(\Pi p_i)_{1 \le i \le n}$ may be uniformly divided into cycles of length $(\Pi p_i)_{1 \le i \le n}/(p_j - 1)$ for each prime $p_j$, each with $f_{n-1}((\Pi(p_i - 1))_{1 \le i \le n-1})$ n-primes.

**Corollary 4.** To that we may add that the number of such cycles is $(p_i - 1)$ for each $p_i$, $i = 1 \cdots n$, resulting in $\sum_{1 \le i \le n}(p_i - 1)$ such predictable and easily calculatable intervals, as illustrated in Table 2.

**Table 2**
**Distribution of 10-Primes**

$\Pi p_i = p_1 \cdot p_2 \cdot p_3 \cdots p_{10} = 2 \cdot 3 \cdot 5 \cdot 7 \cdot 11 \cdot 13 \cdot 17 \cdot 19 \cdot 23 \cdot 29 = 6{,}469{,}693{,}230$

| n | Primes $p_i$ | Number of Intervals $(p_i - 1)$ | Size of Each Interval $(\Pi p_i)_{1 \le i \le n}/(p_j - 1)$ | Number of n-primes per interval 1,021,870,080 / $(p_i - 1)$ |
|---|---|---|---|---|
| 1. | 2 | 1 | 6,469,693,230.0 | 1,021,870,080 |
| 2. | 3 | 2 | 3,234,846,615.0 | 510,935,040 |
| 3. | 5 | 4 | 1,617,423,308.0 | 255,467,520 |
| 4. | 7 | 6 | 1,078,282,205.0 | 170,311,680 |
| 5. | 11 | 10 | 646,969,323.0 | 102,187,008 |
| 6. | 13 | 12 | 539,141,102.5 | 85,155,840 |
| 7. | 17 | 16 | 404,355,826.9 | 63,866,880 |
| 8. | 19 | 18 | 359,427,401.7 | 56,770,560 |
| 9. | 23 | 22 | 294,076,965.0 | 46,448,640 |
| 10. | 29 | 28 | 231,060,472.5 | 36,495,360 |

Here are a couple examples to interpret the table.

- From line 5, we can see that there are 10 intervals, each of size 646,969,230, each one with 102,187,008 10-primes.
- From line 9, we can see that there are 22 intervals, each of size 294,076,965, each one with 46,448,640 10-primes.

**6. A Generalization of N-Primes to Arbitrary Sets of M-Primes**

We can generalize all our results by describing what I will now call M-primes, where M is any set of prime numbers. M-primes would be defined as the integers not divisible by any of the primes in M.

Let's start with the case of M being any finite subset of the primes. All of the concepts and formulas we have been using for the first n primes are generalizable for the M-primes. For example:

- **Theorem 1a.** The distribution of M-primes is repeated every $\Pi p_i$ interval, $p_i \in M$. That is, if x is an M-prime, so is $K \cdot \Pi p_i + x$ for every integer K, $p_i \in M$

- **Theorem 2a.** The M-primes are distributed symmetrically within $P_M$. That is, if x is an M-prime, so is $K \cdot P_M - x$ for every integer K

- **Theorem 3a.** The number of M-primes in $P_M = \Pi(p_i-1)$, $p_i \in M$

- **Corollary 1a.** For any integer $K \geq 0$,
  $f_M(K \cdot \Pi p_i \pm x) = f_M(K \cdot \Pi p_i) \pm f_M(x) = K \cdot \Pi(p_i - 1) \pm f_M(x)$, $p_i \in M$

- **Theorem 4a.** The number of sets of twin M-primes in $P_M = \Pi(p_i-2)$, $p_i \in M$, $(i \neq 2)$. (A twin M-prime is defined to be an integer x such that x–1 and x+1 are both M-primes.)

- **Corollary 2a.** For any integer $a \pm 0 \pmod{P_M}$ and any integer $K \geq 0$ there are $K \cdot \Pi(p_i - 2)$, $p_i \in M$, $(i \neq 2)$ twin M-primes integers $x$ in $K \cdot P_n$ such that $x - 1, x + 1$ is a pair of twin n-primes.

- **Corollary 3a.** For any integer $a \pm 0 \pmod{P_M}$ there are $\Pi(p_i - 2)$, $p_i \in M$, $(p_i \neq 2)$ integers $x$ in $P_M$ such that the pair $x - a, x + a \pmod{P_M}$, are both n-primes.

- **Theorem 5a.** A prime p not $\in M$ divides one and only one n-prime in the series $K \cdot \Pi p_i + a$, $p_i \in M$

- **Theorem 6a.** $f_M((\Pi p_i)/(p_j - 1)) = f_{M-j}(\Pi p_i)$ for any n.
  $\phantom{xx}_{i \in M} \phantom{xx}_{i \in M-j}$

- **Theorem 7a.** There are uniform intervals K of size $P_M/(p_i - 1)$ that each contain $\Pi(p_i - 1)$, $i \in M$, $i \neq j$ for any $j \in M$, M-primes and for any integer K. (We call this the "M-prime distribution formula.")

- **Theorem 8, Meissel's formula generalized.** (Repeated here for ease of reference) $f_M(x) = f_{M-pi}(x) - f_{M-pi}(x/p_i)$ for any $p_i \in M$. (Note, Meissel's formula is expressed in terms of the first n primes. This expression generalizes it arbitrary finite sets of primes.

- **Theorem 9a.** $f_M(K \cdot (\Pi p_i)/(p_j - 1)) = f_{M-pj}(K \cdot (\Pi p_i))$
  $\phantom{xxxxxxx}_{i \in M} \phantom{x}_{i \in M-pj}$

- **Corollary 4a.** The total number of intervals with easily calculable and recurring numbers of M-primes predicted by theorem 9a is $\Sigma(p_i - 1)$, $p_i \in M$.

## 7. M-Primes Generated by Infinite Sets of Primes

Where M is an infinite set of prime numbers (e.g., every second prime, or other calculatable recurring subset, or even a random, but still infinite, set of prime numbers), the theorems above don't make sense because the formulas would result in meaningless infinite products or sums. However, given any number x, the theorems hold up for all primes $\leq x$, and perhaps more interestingly for all primes $\leq$ the square root of x. And in that respect, they hold for any infinite subsets of the whole set of integers.

## 8. A Further Generalization to Relatively Prime Factors of Any Integer or to Any Infinite Sets of Relatively Prime Integers

We may safely observe here the following generalization of the theorems and corollaries of section 6 and reflections on infinite sets of primes in section 7, without meticulously repeating the proofs but simply by adapting the statements and associated proofs to the following generalization—namely, all those theorems and associated corollaries hold if we apply them to *any relatively prime factors of any integer x*, or to *any infinite set of relatively prime integers*. We define relatively prime factors of an integer as any factors of an integer which are not themselves divisible by a common prime factor. For example, 20 and 2,783 are factors of 55,660. Neither of them are primes, but they are relatively prime since they share no prime factors ($20 = 2^2 \cdot 5$, and $2,783 = 11^2 \cdot 23$). Thus, we may use the theorems in section 6 to predict the number of integers $< 55,660$ that are not divisible by 20 or 2,783, being confident as well that such integers symmetrically located between 1 and 55,660, that that pattern will repeat itself with every multiple of 55,600, etc.

## 9. M-Primes and Euler's ϕ-function

N-primes and M-primes build upon and reflect a generalization of Euler's ϕ-function, while at the same time bringing refinements and more penetrating descriptions to it. The fundamental theorem of arithmetic states that every integer x has a unique factorization of primes such that

$$x = p_1^{k_1} p_2^{k_2} p_3^{k_3} \cdots p_n^{k_n}$$

Euler's ϕ function describes the number of integers ≤ an integer x that are relatively prime to it. The formula is written in several equivalent formats. Here we use the following representation:

$$\phi(x) = x \cdot \frac{(p_1 - 1)(p_2 - 1)(p_3 - 1) \cdots (p_n - 1)}{p_1 \cdot p_2 \cdot p_3 \cdots p_n}$$

In this formulation, we do not mean that $p_1, p_2, \ldots p_n$ are the first n primes, but the set M of unique prime dividers of x.

We may write this as

$$\phi(x) = \frac{x}{p_1 \cdot p_2 \cdot p_3 \cdots p_n} \cdot (p_1 - 1)(p_2 - 1)(p_3 - 1) \cdots (p_n - 1),$$

noting that $\frac{x}{p_1 \cdot p_2 \cdot p_3 \cdots p_n}$ is an integer K, since each $p_i$ is a divisor of x.

Therefore,

$$\begin{aligned}\phi(x) &= K \cdot (p_1 - 1)(p_2 - 1)(p_3 - 1) \cdots (p_n - 1), \\ &= K \cdot f_M(p_1 \cdot p_2 \cdot p_3 \cdots p_n) \\ &= f_M(K \cdot p_1 \cdot p_2 \cdot p_3 \cdots p_n) \\ &= f_M(x)\end{aligned}$$

(12) Thus, $\phi(x) = f_M(x)$, where M is the set of prime divisors of x.

On the surface, then, the values of ϕ(x) and $f_M(x)$ are identical. However, ϕ(x) brings with it a simple numerical solution, whereas $f_M(x)$ is "richer" than ϕ(x) in that it provides information not just about the total number of integers that are ≤ x and that are relatively prime to it, but also easily calculatable information about the distribution of the integers that ϕ(x) counts, including all the cycles and waves of integers that are relatively prime to x and that we have discussed throughout this paper—especially the details supplied by the M-prime distribution formulas described above.

**10. Legendre's Function**

Our analysis to this point has enabled us to describe the distribution of n-primes as predictable, recurring sets and subsets of uniformly long cycles and sub cycles along the positive real number line. A more general question might be: Given any set of prime numbers and any real number x, how many n-primes are there that are ≤ x? Or, as we might put it here, is there a general formula for $f_n(x)$ for any real number x? This question was answered by Adrien-Marie Legendre (1752 – 1833) with his formula

(10)  $f_n(x) = [x] - \Sigma [x/p_i] + \Sigma [x/p_i \cdot p_j] - \Sigma [x/p_i \cdot p_j \cdot p_k]$
          $+ \ldots \Sigma [x/p_1 \cdot p_{23} \ldots p_n]$

where in each summation term the primes in the denominator ranges over the first n primes as singles, then pairs, then triplets, etc. The last term will be either positive or negative depending on n (positive if n is an even number, negative if odd). We offer it here without much elaboration simply to round out our presentation of n-primes. We comment only that in its most general applications it requires $2^n$ separate calculations, and does not in and of itself reveal the various patterns of the distribution of n-primes which are the subject of this paper. That said, it is a remarkable and intriguing formula, and one that is more easily used in this day of high speed computers.

Legendre proved his theorem by the general inclusion/exclusion combinatorial theorem. It may also be derived from Meissel theorem (1) as follows:

(11)  The number of integers ≤ x not divisible by a prime $p_1 = [x] - [x/p_1]$. Successive applications of Meissel's theorem ($f_n(x) = f_{n-1}(x) - f_{n-1}(x/p_n)$) then yields Legendre's function.

It is worth noting that Legendre's formula generalizes to any set M of primes, not just the first n primes. Similarly, it holds for all infinite sets of primes, where products of primes in M are irrelevant to the formula once those products exceed x. Of particular relevance here is that this formula used in conjunction with even the most elementary personal computers enables us to easily verify the results of the theorems in this paper.

**11. N-Primes and M-Primes Governed by Cross Product Modular Arithmetic**

The above analysis takes advantage of cross product vector representation of integers, Euler's ϕ function, and the Chinese Remainder Theorem, which together have broad implications for the distribution of and relationships among n-primes. However, the representation of integers as cross product vectors within their first cycle (i.e., within the range $0 \le x \le p_i - 1$ for each i), along with the associated modular arithmetic, greatly facilitates the understanding of relationships of n-primes with the integers of their first cycle and among themselves. Especially see M. Schmitt's thorough treatment of this, focusing on modular structures and calculations for the first n primes[iv].

Of particular interest is that within this Ring, the n-primes themselves form a commutative Abelian subgroup under *cross multiplication* (mod $\Pi p_i$, $1 \leq i \leq n$). This subgroup consists of all vectors ($a_1, a_2, ..., a_n$) such that $a_i \neq 0$ for any $a_i$. For such vectors, that is, for n-primes a, b, and c, with

$$a = (a_1, a_2, ..., a_n)$$
$$b = (b_1, b_2, ..., b_n)$$
$$c = (c_1, c_2, ..., c_n)$$
$$ab = (a_1b_1, a_2b_2, ..., a_nb_n)$$
$$1 = (1, 1, ..., 1)$$

we have
- Closure: a x b = ab
- Associativity: (a x b) x c = a x (b x c)
- An identity element: 1, i.e., a x 1 = a
- Commutativity: a x b = b x a, and
- An inverse element: for any element a there is an element b such that a x b = 1, with the inverse of an element a, sometimes designated 1/a or $a^{-1}$. For each element $a_i$ of the vector ($a_1, a_2, ..., a_n$), there is an element $a_i^{-1}$ such that $a_i \times a_i^{-1} = m_i p_{i\,+}\, 1$ for some integer $m_i$.

It is worth clarifying that, depending on the topic under discussion, we have sometimes switched between the modular arithmetic described here and normal integer or rational number arithmetic, depending on the topic under discussion. For example, theorem's 2 – 4 above can be interpreted either way. However, we sometimes needed to analyze n-primes as integers embedded in the larger domains of real numbers, particularly when calculating the dividing lines between ranges of the subsets of n-primes that we discuss. Hopefully, the context made it clear which domain of reference was at play.

**12. Conclusion**

More than 2300 years ago, without the benefit of high speed computers, Eratosthenes of Cyrene peered with his mind's eye far down the endless number line, satisfying some of his own curiosity about prime numbers, while teasing some 600 generations of aspiring mathematicians to learn more about them. This paper does not reveal much more about the primes, but hopefully sheds some light on the survivors of the infinitely destructive rampages unleashed by his sieve. In their path, and even with our elementary methods, we find well defined, predictably calculatable, symmetrical cycles and sub-cycles of survivors. Such patterns and their related formulas will not be extinguished any more than the primes themselves, no matter how many times Eratosthenes's sieve is applied.

georgefgrob@cs.com, MA, Georgetown University, 1969
mschmitt@db12.de, diploma mathematician graduated in Mathematical Institute, University of Cologne, 1997